\documentclass[a4paper,12pt]{article}

\usepackage{color,makeidx,bm,latexsym,amscd,amsmath,amssymb,stmaryrd,amsthm,float,graphicx,psfrag}

\theoremstyle{plain}
\newtheorem{theorem}{Theorem}[section]
\newtheorem{corollary}[theorem]{Corollary}
\newtheorem{lemma}[theorem]{Lemma}
\newtheorem{conjecture}[theorem]{Conjecture}

\theoremstyle{definition}
\newtheorem{definition}[theorem]{Definition}
\newtheorem{remark}[theorem]{Remark}
\newtheorem{example}[theorem]{Example}
\newtheorem{proposition}[theorem]{Proposition}
\newtheorem{assumption}[theorem]{Assumption}

\newcommand{\bF}{\ensuremath{\mathbb{F}}}

\newcommand{\bK}{\ensuremath{\mathbb{K}}}

\newcommand{\bZ}{\ensuremath{\mathbb{Z}}}

\newcommand{\scL}{\ensuremath{\mathcal{L}}}

\newcommand{\A}{\mathcal{A}}
\newcommand{\C}{\mathbb{C}}

\newcommand{\Q}{\mathbb{Q}}
\newcommand{\R}{\mathbb{R}}
\newcommand{\Z}{\mathbb{Z}}

\newcommand{\asc}{\operatorname{asc}}

\newcommand{\geneul}{R}
\newcommand{\ehr}{L}
\newcommand{\ehrser}{\operatorname{Ehr}}
\newcommand{\alcov}{{\sigma_\Phi^\circ}}
\newcommand{\baralcov}{\overline{\alcov}}

\newcommand{\cwtl}{Z}
\newcommand{\crtl}{\check{Q}}

\newcommand{\quasi}{\operatorname{quasi}}

\newcommand{\rad}{\operatorname{rad}}

\renewcommand{\Re}{\operatorname{Re}}

\newcommand{\barH}{\overline{H}}
\newcommand{\Zq}{\Z/q\Z}

\title{Characteristic polynomials of Linial arrangements 
for exceptional root systems}
\author{Masahiko Yoshinaga\thanks{Department of Mathematics, 
Hokkaido University, North 10, West 8, Kita-ku, 
Sapporo 060-0810, JAPAN 
E-mail: yoshinaga@math.sci.hokudai.ac.jp}}
\date{\today}
\begin{document}
\maketitle

\begin{abstract} 
The (extended) Linial arrangement $\mathcal{L}_{\Phi}^m$ is a certain finite 
truncation of the affine Weyl arrangement of a root system 
$\Phi$ with a parameter $m$. 
Postnikov and Stanley conjectured that all roots of 
the characteristic polynomial of $\mathcal{L}_{\Phi}^m$ have the 
same real part, and this has been proved for the root systems of 
classical types. 

In this paper we prove that the conjecture is true for exceptional root 
systems when the parameter $m$ is sufficiently large. 

The proof is based on representations of 
the characteristic quasi-polynomials 
in terms of Eulerian polynomials. 


\end{abstract}

\tableofcontents

\section{Introduction}
\label{sec:intro}

\subsection{Background}
\label{subsec:background}

A hyperplane arrangement $\A=\{H_1, \dots, H_n\}$ is a finite 
collection of affine hyperplanes in an $\ell$-dimensional vector space 
$\bK^\ell$. Despite its simplicity, the theory of hyperplane 
arrangements has fruitful connections with many areas in mathematics 
(\cite{ot, st-lect}). 
One of the most important invariant of an arrangement $\A$ 
is the 
\emph{characteristic polynomial} $\chi(\A, t)\in\Z[t]$. 
Indeed the characteristic polynomial is related to several 
other invariants, such as 
the Poincar\'e polynomial of the complexified complement 
$M(\A)$ \cite{os}, the number of chambers for real arrangements 
\cite{zas-face}, the number of $\bF_q$-rational points 
\cite{cra-rot, ter-jac}, Chern classes of certain vector bundles 
\cite{mus-sch, alu}, and lattice points countings 
\cite{bl-sa, ktt-cent, ktt-noncent, ktt-quasi, yos-worp}. 

\subsection{Main results}
\label{subsec:result}

Let $V=\R^\ell$ be an $\ell$-dimensional Euclidean space. 
Let $\Phi\subset V^*$ be an irreducible root system. 
Fix a positive system $\Phi^+\subset\Phi$. For a positive root 
$\alpha\in\Phi^+$ and $k\in\Z$, define 
\[
H_{\alpha, k}=\{x\in V\mid \alpha(x)=k\}. 
\]
The set of all such hyperplanes is called the affine Weyl arrangement. 
Finite truncations of the affine Weyl arrangement have received 
considerable attention (\cite{ath-adv, ath-survey, ath-lin, ath-gen, ede-rei, 
ps-def, shi-kl, ter-multi, yos-char}). Among others, 
the (extended) Linial arrangement $\scL_\Phi^m$ is defined by 
\[
\scL_\Phi^m=\{H_{\alpha, k}\mid
\alpha\in\Phi^+, k=1, 2, \dots, m\}, 
\]
(where $\scL_\Phi^0=\emptyset$ by convention). 
In \cite{ps-def}, Postnikov and Stanley studied combinatorial 
aspects of Linial arrangements. 
They posed the following conjecture. 

\begin{conjecture}
\label{conj:rh}
(\cite[Conjecture 9.14]{ps-def}) 
Suppose $m\geq 1$. Then every root $\alpha\in\C$ of the equation 
$\chi(\scL_{\Phi}^{m}, t)=0$ satisfies 
$\Re \alpha=\frac{mh}{2}$, where $h$ denotes the Coxeter number of $\Phi$. 
\end{conjecture}
The conjecture was verified for $\Phi=A_\ell$ 
by Postnikov and Stanley \cite{ps-def}, 
and for $\Phi=B_\ell, C_\ell$, and $D_\ell$ by 
Athanasiadis (\cite{ath-lin}). These works are based on 
explicit representations of $\chi(\scL_{\Phi}^{m}, t)$ for the corresponding 
root systems. (The case $\Phi=G_2$ is also easy.) 

For exceptional root 
systems, some partial answers have recently been reported in \cite{yos-worp}. 
Namely, for $\Phi\in\{E_6, E_7, E_8, F_4\}$, 
Conjecture \ref{conj:rh} has been verified when the parameter $m>0$ satisfies 
\[
m\equiv -1
\left\{
\begin{array}{ll}
\mod 6, & \Phi=E_6, E_7, F_4, \\
\mod 30, & \Phi=E_8. 
\end{array}
\right.
\]

The purpose of this paper is to prove Conjecture \ref{conj:rh} 
for exceptional root systems 
when $m\gg 0$. The main result is the following. 

\begin{theorem}
\label{thm:main}
(Corollary \ref{cor:main}) 
Let $\Phi\in\{E_6, E_7, E_8, F_4\}$. 
Suppose $m\gg 1$. 
Then, every root $\alpha\in\C$ of the equation 
$\chi(\scL_{\Phi}^{m}, t)=0$ satisfies 
$\Re \alpha=\frac{mh}{2}$. 
\end{theorem}

\subsection{What makes roots lie on a line?}

The proof of Theorem \ref{thm:main} relies on the 
expression of the characteristic quasi-polynomial 
$\chi_{\quasi}(\scL_{\Phi}^m, q)$ in terms of 
the Ehrhart quasi-polynomials and Eulerian polynomials 
developed in \cite{yos-worp}. (See \S \ref{sec:pre}.) 
However, the key result that enables us to 
conclude ``having the same real part'' is the following elementary lemma. 

\begin{lemma}
\label{lem:elem}
Let $f(t)\in\R[t]$. Suppose $M$ is a real number that satisfies 
the inequality
\[
M>2\cdot\max\{\Re\alpha\mid\alpha\in\C, f(\alpha)=0\}. 
\]
Let $\omega\in\C$ be a complex number with $|\omega|=1$. 
Then, the root $\alpha$ of the equation 
\[
f(t)-\omega\cdot f(M-t)=0
\]
satisfies $\Re \alpha=\frac{M}{2}$. 
\end{lemma}
\begin{proof}
Set $f(t)=a(t-\alpha_1)(t-\alpha_2)\cdots(t-\alpha_n)$, 
($a\neq 0$). As $f(t)$ is a real polynomial, $\overline{\alpha_i}$ is 
also a root of $f(t)$. Set $\beta_i=M-\overline{\alpha_i}$. 
Then, $\alpha_i$ and $\beta_i$ are symmetric with respect to the line 
$\Re =\frac{M}{2}$, and we have $f(M-t)=(-1)^n\cdot\prod_{i=1}^n(t-\beta_i)$. 
If $\Re z<\frac{M}{2}$, then $|z-\alpha_i|<|z-\beta_i|$ for all $i$, 
and hence $|f(z)|<|f(M-z)|$. 
Similarly, $\Re z>\frac{M}{2}$ implies $|f(z)|>|f(M-z)|$. 
Therefore, $f(z)=\omega f(M-z)$ implies that 
$\Re z=\frac{M}{2}$. 
\end{proof}

The basic strategy of the proof of Theorem \ref{thm:main} is to construct 
$F^{(m)}(t)\in\Q[t]$ such that 
\[
\chi(\scL_{\Phi}^m, t)=F^{(m)}(t)+(-1)^\ell\cdot F^{(m)}(mh-t), 
\]
where $\ell$ is the rank of $\Phi$, 
then apply Lemma \ref{lem:elem}. 

\medskip

The remainder of this paper is organized as follows. 
In \S \ref{sec:pre}, 
we recall the notion of the characteristic quasi-polynomial 
$\chi_{\quasi}(\A, q)$ for an integral arrangement $\A$. 
The characteristic quasi-polynomial $\chi_{\quasi}(\scL_\Phi^m, q)$ 
of the Linial arrangement $\scL_\Phi^m$ can be expressed 
in terms of the 
Ehrhart quasi-polynomial $L_\Phi(t)$ of the fundamental alcove 
and the Eulerian polynomial $R_\Phi(t)$ of $\Phi$. 
The relation between these objects is described as follows 
(Theorem \ref{thm:yos-worp}) 
\[
\chi_{\quasi}(\scL_\Phi^m, q)=R_{\Phi}(S^{m+1})L_{\Phi}(q), 
\]
where $S$ is the shift operator. 
In \S \ref{sec:limpoly}, 
using the symmetry of the Eulerian polynomial 
(Proposition \ref{prop:properties}), we introduce the truncated 
characteristic quasi-polynomial 
$\chi_{\quasi}^{1/2}(\scL_\Phi^m, q)$, which satisfies 
\[
\chi_{\quasi}(\scL_{\Phi}^m, q)=\chi_{\quasi}^{1/2}(\scL_{\Phi}^m, q)+
(-1)^\ell\chi_{\quasi}^{1/2}(\scL_{\Phi}^m, mh-q), 
\] 
(Proposition \ref{prop:decquasi}). Since these functions are 
quasi-polynomials, it does not make sense to consider the roots. 
However, we will see that the limit 
\[
F_\Phi(q):=R_\Phi^{1/2}(S)q^\ell=
\lim_{m\rightarrow\infty}
\frac{\chi_{\quasi}^{1/2}(\scL_\Phi^m, mq)}{m^\ell}
\]
becomes a polynomial in $q$ (Proposition \ref{prop:normalim}). 
The location of the zeros of $F_\Phi(t)$ is crucial for $m\gg 0$. 
Indeed, we will check, case-by-case, that the real parts of 
the zeros are less than $\frac{h}{2}$ (Proposition \ref{prop:Remax}). 

Because the proof for the quasi-polynomials is complicated, we will first 
give a simplified ``polynomial  version'' of 
the main result as a ``toy-case'' in \S \ref{sec:toy}. 

In \S \ref{subsec:settings}, we summarize those properties of 
quasi-polynomials and Eulerian polynomials that are necessary 
for the proof of the main result. This is mainly to simplify 
the notation. In \S \ref{subsec:asympt}, we present a weaker 
version of the main result for the asymptotic behavior of the real parts of 
roots. In \S \ref{subsec:exact}, we prove the main result. 

\section{Preliminaries}

\label{sec:pre}

\subsection{Quasi-polynomials with the $GCD$-property} 

\label{subsec:quasi}

A map $F:\Z\longrightarrow\C$ is called a \emph{quasi-polynomial} 
with a period $\rho>0$ if $F(q)$ can be expressed as a polynomial in $q$ 
that depends only on the residue class $q \mod \rho$. 
In other words, there exist 
polynomials $f_0, f_1, \dots, f_{\rho-1}\in\C[t]$ such that 
\[
F(q)=f_i(q) 
\]
if $q\equiv i\mod\rho$. 
The polynomials $f_0, \dots, f_{\rho-1}$ are called the 
constituents of $F$. The period $\rho$ is said to be the minimal period if $F$ 
does not have smaller periods than $\rho$. The quasi-polynomial $F$ can be 
expressed as 
\begin{equation}
\label{eq:periodicpoly}
F(q)=c_0(q)q^d+c_1(q)q^{d-1}+\dots+c_d(q), 
\end{equation}
where $c_i:\Z\longrightarrow\C$ is a periodic function 
with a period $\rho$ (i.e., $c_i(q+\rho)=c_i(q)$ for all $q\in\Z$). 

We say that the quasi-polynomial $F$ has a constant leading term if 
$c_0(q)$ in (\ref{eq:periodicpoly}) is a nonzero constant function. 
In this case, $d$ is called the degree of the quasi-polynomial $F$. 

The quasi-polynomial $F$ is said to have the $GCD$-property 
if the constituents satisfy: 
$\gcd(i, \rho)=\gcd(j, \rho)\Longrightarrow 
f_i(t)=f_j(t)$. 

\begin{remark}
In this paper, we distinguish the roles of the variables $q$ and $t$. 
The variable $q$ always runs through $\Z$ (or $\Z_{>0}$), 
whereas $t$ runs through $\R$ or $\C$. 
Under this convention, the variable of 
a quasi-polynomial should be $q$, and its constituents may have 
a variable $t$. 
\end{remark}

\subsection{Characteristic quasi-polynomials}

\label{subsec:charquasipoly}

Let $\A=\{H_1, \dots, H_n\}$ be an arrangement of affine hyperplanes 
in $\R^n$. Throughout this paper, we assume that the hyperplanes are 
defined over $\Z$. More precisely, there exists an integral linear equation 
\[
\alpha_i(x_1, \dots, x_\ell)=a_{i1}x_1+\dots+a_{i\ell}x_\ell+b_i 
\]
($a_{ij}, b_i\in\Z$) that satisfies $H_i=\alpha_i^{-1}(0)\subset\R^\ell$. 
For an arrangement $\A$, 
we can associate the modulo $q>0$ complement: 
\[
M_q(\A)=\left(\Zq\right)^\ell\smallsetminus
\bigcup_{i=1}^n\barH_i, 
\]
where $\barH_i=\{x\in\left(\Zq\right)^\ell\mid\alpha_i(x)\equiv 0\mod q\}$. 

The following theorem was given by Kamiya, Takemura and Terao. 

\begin{theorem}
(\cite{ktt-cent, ktt-noncent, ktt-quasi}) 
$\#M_q(\A)$ is a quasi-polynomial with the $GCD$-property for 
sufficiently large $q\gg 0$. 
\end{theorem}
We denote the quasi-polynomial by 
$\chi_{\quasi}(\A, q)$, which is called the \emph{characteristic 
quasi-polynomial} of $\A$. 
The characteristic quasi-polynomial 
has a constant leading term; hence, it is of the form 
\[
\chi_{\quasi}(\A, q)=q^\ell+c_1(q)\cdot q^{\ell-1}+\dots+c_\ell(q), 
\]
where $c_i:\Z\longrightarrow\Z$, $i=1, \dots, \ell$ are periodic functions. 
It is also known that the prime constituent 
of $\chi_{\quasi}(\A, q)$ is equal to the characteristic polynomial of $\A$ 
\cite{ath-adv, ath-lin}, i. e., 
the characteristic polynomial $\chi(\A, t)$ has the form 
\[
\chi(\A, t)=t^\ell+c_1(1)\cdot t^{\ell-1}+\dots+c_\ell(1). 
\]

\subsection{Eulerian polynomials for root systems}

\label{subsec:epoly}

We first recall the terminology of \cite{bour, hum}. 

Let $V=\R^\ell$ be the Euclidean space with inner product 
$(\cdot, \cdot)$. Let $\Phi\subset V$ be an irreducible 
root system with exponents $e_1, \dots, e_\ell$, 
Coxeter number $h$, and Weyl group $W$. 
For any integer $k\in\Z$ and $\alpha\in\Phi^+$, 
the affine hyperplane $H_{\alpha, k}$ is defined by 
\begin{equation}
\label{eq:affinehyperp}
H_{\alpha, k}=\{x\in V\mid 
(\alpha, x)=k\}. 
\end{equation}

Fix a positive system $\Phi^+\subset \Phi$ and 
the set of simple roots $\Delta=\{\alpha_1, \dots, \alpha_\ell\}
\subset\Phi^+$. The highest root, denoted by 
$\widetilde{\alpha}\in\Phi^+$, can be expressed as a linear 
combination $\widetilde{\alpha}=\sum_{i=1}^\ell c_i\alpha_i$ 
($c_i\in\Z_{>0}$). We also set $\alpha_0:=-\widetilde{\alpha}$ and 
$c_0:=1$. Then, we have the linear relation 
\begin{equation}
\label{eq:linrel}
c_0\alpha_0+c_1\alpha+\dots+c_\ell\alpha_\ell=0. 
\end{equation}
The coweight lattice $\cwtl(\Phi)$ and the coroot lattice $\crtl(\Phi)$ 
are defined as 
\[
\begin{split}
\cwtl(\Phi)
&=
\{x\in V\mid (\alpha_i, x)\in\Z, \alpha_i\in\Delta\}, \\
\crtl(\Phi)
&=
\sum_{\alpha\in\Phi}\Z\cdot
\frac{2\alpha}{(\alpha, \alpha)}. 
\end{split}
\]
The coroot lattice $\crtl(\Phi)$ is a finite index subgroup 
of the coweight lattice $\cwtl(\Phi)$. The index 
$\#\frac{\cwtl(\Phi)}{\crtl(\Phi)}=f$ is called the 
\emph{index of connection}. 

Let $\varpi_i^\lor\in\cwtl(\Phi)$ be the dual basis of the 
simple roots $\alpha_1, \dots, \alpha_\ell$, that is, 
$(\alpha_i, \varpi_j^\lor)=\delta_{ij}$. Then, $\cwtl(\Phi)$ is 
a free abelian group generated by 
$\varpi_1^\lor, \dots, \varpi_\ell^\lor$. We also have 
$c_i=(\varpi_i^\lor, \widetilde{\alpha})$. 

Each connected component of $V\smallsetminus 
\bigcup\limits_{\substack{\alpha\in\Phi^+\\ k\in\Z}}H_{\alpha, k}$ 
is an open simplex, called an \emph{alcove}. 
Define the fundamental alcove $\alcov$ by 
\[
\alcov
=
\left\{
x\in V
\left|
\begin{array}{ll}
(\alpha_i, x)>0,&(1\leq i\leq \ell)\\
(\widetilde{\alpha}, x)<1&
\end{array}
\right.
\right\}
\]
The closure 
$\baralcov=\{x\in V\mid (\alpha_i, x)\geq 0\ 
(1\leq i\leq \ell),\ (\widetilde{\alpha}, x)\leq 1\}$ 
is a simplex with vertices $0, \frac{\varpi_1^\lor}{c_1}, \dots, 
\frac{\varpi_\ell^\lor}{c_\ell}\in V$. 
The supporting hyperplanes of facets of $\baralcov$ are 
$H_{\alpha_1, 0}, \dots, H_{\alpha_\ell, 0}, H_{\widetilde{\alpha}, 1}$. 

Using the linear relation (\ref{eq:linrel}), 
we define the function $\asc:W\longrightarrow\bZ$. 

\begin{definition}
\label{def:asc}
Let $w\in W$. Then, $\asc(w)$ is defined by 
\[
\asc(w)=\sum_{\substack{0\leq i\leq \ell\\ w(\alpha_i)>0}}c_i. 
\]
\end{definition}

\begin{definition}
\label{def:geneul}
The generalized Eulerian polynomial $R_\Phi(x)$ is defined by 
\[
R_\Phi(x)=\frac{1}{f}\sum_{w\in W}x^{\asc(w)}. 
\]
\end{definition}
The following proposition gives some basic properties of 
$\geneul_\Phi(x)$. 
\begin{proposition}
\label{prop:properties}
(\cite{lp-alc2})
\begin{itemize}
\item[(1)] $\deg \geneul_\Phi(x)=h-1$. 
\item[(2)] (Duality) $x^h\cdot\geneul_\Phi(\frac{1}{x})=\geneul_\Phi(x)$. 
\item[(3)] $\geneul_\Phi(x)\in\bZ[x]$. 
\item[(4)] $\geneul_{A_\ell}(x)$ is equal to the classical 
Eulerian polynomial. (See \cite{comtet, foa-hist, st-ec1} for 
classical Eulerian polynomials.) 
\end{itemize}
\end{proposition}

The polynomial $\geneul_\Phi(x)$ was introduced by 
Lam and Postnikov in \cite{lp-alc2}. 
They proved 
that $\geneul_\Phi(x)$ can be expressed in terms of cyclotomic 
polynomials and classical Eulerian polynomials. 
\begin{theorem}
\label{thm:lp}
(\cite[Theorem 10.1]{lp-alc2}) Let $\Phi$ be a root system of rank $\ell$. 
Then, 
\begin{equation}
\label{eq:lp}
\geneul_\Phi(x)=[c_0]_x\cdot [c_1]_x\cdot [c_2]_x\cdots [c_\ell]_x\cdot
\geneul_{A_\ell}(x), 
\end{equation}
where $[c]_x=\frac{x^c-1}{x-1}$. 
\end{theorem}
We will give an alternative proof of Theorem \ref{thm:lp} in 
\S \ref{subsec:cqpLinial} using Ehrhart series.

\subsection{Ehrhart quasi-polynomials for root systems}

\label{subsec:eqpoly}

It is known that the number of lattice points 
$L_\Phi(q)=\#\left(q\cdot\baralcov\cap Z(\Phi)\right)$ in the delate 
$q\cdot\baralcov$ is a quasi-polynomial in $q$, 
called the \emph{Ehrhart quasi-polynomial} of $\baralcov$. 
(See \cite{be-ro, st-ec1} for details on Ehrhart theory.) 
Suter \cite{sut} 
explicitly computed the 
Ehrhart quasi-polynomial $\ehr_{\baralcov}(q)$. 
Several useful conclusions may be summarized as follows. 
\begin{theorem}
\label{thm:suter}
(Suter \cite{sut}) 
\begin{itemize}
\item[(i)] 
The Ehrhart quasi-polynomial $\ehr_{\baralcov}(q)$ 
has the $GCD$-property. 
\item[(ii)] 
$\ehr_{\baralcov}(q)$ has a constant leading term whose 
leading coefficient is $\frac{f}{|W|}$. 
\item[(iii)] 
The minimal period is 
$\widetilde{n}=\gcd(c_1, c_2, \dots, c_\ell)$. 
(See Table \ref{fig:table} for explicit values.) 
\item[(iv)] 
If $q\in\Z$ is relatively prime to the period $\widetilde{n}$, then 
\[
\ehr_{\baralcov}(q)=\frac{f}{|W|}(q+e_1)(q+e_2)\cdots(q+e_\ell). 
\]
\item[(v)] $\rad(\widetilde{n})|h$, 
where $\rad(\widetilde{n})=
\prod_{p: \mbox{\scriptsize prime}, p|\widetilde{n}}p$ 
is the radical of $\widetilde{n}$. 
\item[(vi)] The Ehrhart series $\ehrser_{\Phi}(z)$ of $L_\Phi(q)$ is 
\begin{equation}
\label{eq:ehrser}
\ehrser_{\Phi}(z):=\sum_{q=0}^\infty L_\Phi(q)z^q=
\frac{1}{(1-z^{c_0})(1-z^{c_1})\cdots(1-z^{c_\ell})}. 
\end{equation}
\end{itemize}
\end{theorem}

\begin{table}[htbp]
\centering
{\footnotesize 
\begin{tabular}{c|c|c|c|c|c|c|c}
$\Phi$&$e_1, \dots, e_\ell$&$c_1, \dots, c_\ell$&$h$&$f$&$|W|$&$\widetilde{n}$&$\rad(\widetilde{n})$\\
\hline\hline
$A_\ell$&$1,2,\dots,\ell$&$1,1,\dots,1$&$\ell+1$&$\ell+1$&$(\ell+1)!$&$1$&$1$\\
$B_\ell, C_\ell$&$1,3,5,\dots,2\ell-1$&$1,2,2,\dots,2$&$2\ell$&$2$&$2^\ell\cdot \ell!$&$2$&$2$\\
$D_\ell$&$1,3,5,\dots,2\ell-3,\ell-1$&$1,1,1,2,\dots,2$&$2\ell-2$&$4$&$2^{\ell-1}\cdot\ell!$&$2$&$2$\\
$E_6$&$1,4,5,7,8,11$&$1,1,2,2,2,3$&$12$&$3$&$2^7\cdot 3^4\cdot 5$&$6$&$6$\\
$E_7$&$1,5,7,9,11,13,17$&$1,2,2,2,3,3,4$&$18$&$2$&$2^{10}\cdot 3^4\cdot 5\cdot 7$&$12$&$6$\\
$E_8$&$1,7,11,13,17,19,23,29$&$2,2,3,3,4,4,5,6$&$30$&$1$&$2^{14}\cdot 3^5\cdot 5^2\cdot 7$&$60$&$30$\\
$F_4$&$1,5,7,11$&$2,2,3,4$&$12$&$1$&$2^7\cdot 3^2$&$12$&$6$\\
$G_2$&$1,5$&$2,3$&$6$&$1$&$2^2\cdot 3$&$6$&$6$
\end{tabular}
}
\caption{Table of root systems.}
\label{fig:table}
\end{table}

The following proposition is a conclusion of 
the Ehrhart--Macdonald reciprocity. 
\begin{proposition}
\label{prop:hshift}
\begin{equation}
\ehr_{\Phi}(-q)=(-1)^\ell\cdot \ehr_\Phi(q-h). 
\end{equation}
(See \cite[Corollary 3.4]{yos-worp} for a more general formula.) 
\end{proposition}

\subsection{Characteristic quasi-polynomials of Linial arrangements}

\label{subsec:cqpLinial}

Let $S$ be the shift operator that replaces the variable $q$ by 
$q-1$ (or $t$ by $t-1$). 
More generally, the polynomial $\sigma(S)=a_0+a_1S+\cdots+a_dS^d$ 
acts on a function $f(q)$ as 
\[
\sigma(S)f(q)=a_0f(q)+a_1f(q-1)+\cdots+a_df(q-d). 
\]
The characteristic quasi-polynomial $\chi_{\quasi}(\scL_{\Phi}^n, q)$ 
can be expressed in terms of the Eulerian polynomial $\geneul_\Phi(x)$, 
and the Ehrhart quasi-polynomial $\ehr_\Phi(q)$ of the fundamental alcove. 

\begin{theorem}
\label{thm:yos-worp}
(\cite{yos-worp}) 
\begin{equation}
\label{eq:yos-worp}
\chi_{\quasi}(\scL_\Phi^m, q)=
\geneul_\Phi(S^{m+1})L_\Phi(q). 
\end{equation}
\end{theorem}

\begin{remark}
\label{rem:worp}
The formula (\ref{eq:yos-worp}) holds for 
$m=0$ if we consider $\scL_\Phi^0$ to be the empty arrangement. 
In that case, we have 
\begin{equation}
\label{eq:L0}
q^\ell=\geneul_\Phi(S)L_\Phi(q). 
\end{equation}
This can be considered as a root system generalization of 
the so-called Worpitzky identity \cite{wor}. 

Note that Theorem \ref{thm:lp} by Lam and Postnikov follows 
from the Worpitzky identity (\ref{eq:L0}) and Theorem \ref{thm:suter} (vi). 
Indeed, equation (\ref{eq:L0}) is equivalent to 
\begin{equation}
\label{eq:eulerianser}
\geneul_\Phi(z)\ehrser_{\Phi}(z)=\sum_{k=0}^\infty k^\ell z^k. 
\end{equation}
Note that the right-hand side depends only on the rank $\ell$. 
A comparison of (\ref{eq:eulerianser}) with $\Phi=A_\ell$ 
shows that 
\[
\frac{\geneul_\Phi(z)}{(1-z^{c_0})(1-z^{c_1})\cdots(1-z^{c_\ell})}
=
\frac{\geneul_{A_\ell}(z)}{(1-z)^{\ell+1}}. 
\]
This yields Theorem \ref{thm:lp}. 
\end{remark}

\begin{example}
\label{ex:g2}
Let $\Phi=G_2$. 
Since $\geneul_{A_2}(x)=x+x^2$, we have 
\[
\geneul_{G_2}(x)=(1+x)(1+x+x^2)(x+x^2)=x+3x^2+4x^3+3x^4+x^5. 
\]
The closed fundamental alcove $\baralcov$ is 
the convex hull of $0, \frac{\varpi_1^\lor}{2}, 
\frac{\varpi_1^\lor}{3}$. 
The period is $\widetilde{n}=6$. 
The Ehrhart quasi-polynomial is 
\[
\ehr_{\baralcov}(q)=
\left\{
\begin{array}{ll}
\frac{1}{12}(q+1)(q+5), 
&\mbox{ if $q\equiv 1, 5 \mod 6$},\\
&\\
\frac{1}{12}(q+2)(q+4), 
&\mbox{ if $q\equiv 2, 4 \mod 6$},\\
&\\
\frac{1}{12}(q+3)^2, 
&\mbox{ if $q\equiv 3 \mod 6$},\\
&\\
\frac{1}{12}(q^2+6q+12), 
&\mbox{ if $q\equiv 0 \mod 6$.}
\end{array}
\right.
\]
It is known 
that the characteristic quasi-polynomial of the 
Weyl arrangement $\A_{\Phi}=\{H_{\alpha, 0}\mid \alpha\in
\Phi^+\}$ can be expressed as 
$\chi_{\quasi}(\A_\Phi, q)=
(-1)^\ell\frac{\#W}{f}\ehr_{\baralcov}(-q)$ \cite{ath-gen, ktt-quasi}. 
Thus, we have 
\[
\chi_{\quasi}(\A_{G_2}, q)=
\left\{
\begin{array}{ll}
(q-1)(q-5), 
&\mbox{ if $t\equiv 1, 5 \mod 6$},\\
&\\
(q-2)(q-4), 
&\mbox{ if $t\equiv 2, 4 \mod 6$},\\
&\\
(q-3)^2, 
&\mbox{ if $t\equiv 3 \mod 6$},\\
&\\
q^2-6q+12, 
&\mbox{ if $t\equiv 0 \mod 6$.}
\end{array}
\right.
\]

The characteristic quasi-polynomial of the Linial arrangement is 
\[
\chi_{\quasi}(\scL_{G_2}^1, q)=
\left\{
\begin{array}{ll}
q^2-6q+11 & q\equiv 1\mod 2, \\
q^2-6q+14 & q\equiv 0\mod 2. 
\end{array}
\right. 
\]
\end{example}

\section{Limit Polynomials}

\label{sec:limpoly}

\subsection{Normalized limit polynomials}

\label{subsec:normaliz}

Let $f(S)\in\C[S]$ be a polynomial of the shift operator $S$, 
and $g(t)\in\C[t]$. Assume $\deg g(t)=\ell$. 
Let us consider the polynomial 
\begin{equation}
g_m(t):=f(S^{m+1})g(t) 
\end{equation}
for $m\geq 0$. 

\begin{proposition}
\label{prop:normalim}
\begin{equation}
\lim_{m\to\infty}\frac{g_m(mt)}{m^\ell}=f(S)t^\ell. 
\end{equation}
\end{proposition}

\begin{proof}
Write $f(S)=\sum_{k=0}^N a_kS^k$ and 
$g(t)=\sum_{i=0}^\ell c_it^{\ell-i}$ ($c_0\neq 0)$. 
Then, we have, 
\[
g_m(t)=f(S^{m+1})g(t)=
\sum_{k=0}^N\sum_{i=0}^\ell a_k c_i\cdot (t-k(m+1))^{\ell-i}. 
\]
Hence, 
\[
\begin{split}
\lim_{m\to\infty}
\frac{g_m(mt)}{m^\ell}
&=
\lim_{m\to\infty}
\frac{1}{m^\ell}
\sum_{i=0}^\ell\sum_{k=0}^N a_k c_i\cdot (mt-k(m+1))^{\ell-i}. \\
&=
\lim_{m\to\infty}
\sum_{i=0}^\ell\frac{1}{m^i}
\sum_{k=0}^N a_k c_i\cdot (t-k\frac{m+1}{m})^{\ell-i}. \\
&=
\sum_{k=0}^N a_k c_0\cdot (t-k)^{\ell}\\
&=
f(S)t^\ell. 
\end{split}
\]
\end{proof}

\subsection{Truncated Eulerian polynomials}

\label{subsec:trunc}

Suppose $\geneul_\Phi(x)=\sum_{i=1}^{h-1}a_ix^i$. 
Define the truncated Eulerian polynomial $\geneul_\Phi^{1/2}(t)$ by 
\begin{equation}
\geneul_\Phi^{1/2}(x)=
\left\{
\begin{array}{ll}
\sum\limits_{1\leq i<\frac{h}{2}}a_ix^i, & \mbox{ if $h$ is odd},\\
&\\
\sum\limits_{1\leq i<\frac{h}{2}}a_ix^i+\frac{a_{h/2}}{2}x^{h/2}, & \mbox{ if $h$ is even}. 
\end{array}
\right.
\end{equation}

\begin{example}
Let $\Phi=G_2$. Then, 
$\geneul_\Phi^{1/2}(x)=x+3x^2+2x^3$. 
\end{example}

The following is straightforward from 
Proposition \ref{prop:properties} (2). 

\begin{proposition}
\label{prop:dectrunc}
$\geneul_\Phi(x)=\geneul_\Phi^{1/2}(x)+x^h\cdot\geneul_\Phi^{1/2}(x^{-1})$. 
\end{proposition}
Using the truncated Eulerian polynomial $\geneul_\Phi^{1/2}(x)$, we 
define the half characteristic quasi-polynomial 
$\chi_{\quasi}^{1/2}(\scL_\Phi^m, q)$ as follows. 
\begin{equation}
\chi_{\quasi}^{1/2}(\scL_\Phi^m, q):=\geneul_\Phi^{1/2}(S^{m+1})\ehr_\Phi(q). 
\end{equation}
Note that $\chi_{\quasi}^{1/2}(\scL_\Phi^m, q)$ is a quasi-polynomial 
of the period $\widetilde{n}$, however, it does not have the $GCD$-property. 

\begin{example}
For $\Phi=G_2$, we have 
\[
R_{G_2}^{1/2}(S)L_{G_2}(q)=
\left\{
\begin{array}{ll}
\frac{6q^2+10q}{12}&q\equiv 0\mod 3, \\
\frac{6q^2+10q-4}{12}&q\equiv 1\mod 3, \\
\frac{6q^2+10q+4}{12}&q\equiv 2\mod 3, 
\end{array}
\right.
\]
and 
\[
\chi_{\quasi}^{1/2}(\scL_{G_2}^1, q)=
R_{G_2}^{1/2}(S^2)L_{G_2}(q)=
\left\{
\begin{array}{ll}
\frac{3q^2-8q+12}{6} & q\equiv 0\mod 6, \\
\frac{3q^2-8q+5}{6} & q\equiv 1\mod 6, \\
\frac{3q^2-8q+10}{6} & q\equiv 2\mod 6, \\
\frac{3q^2-8q+3}{6} & q\equiv 3\mod 6, \\
\frac{3q^2-8q+14}{6} & q\equiv 4\mod 6, \\
\frac{3q^2-8q+1}{6} & q\equiv 5\mod 6. 
\end{array}
\right.
\]
\end{example}

\begin{proposition}
\label{prop:decquasi}
The half characteristic quasi-polynomial 
$\chi_{\quasi}^{1/2}(\scL_{\Phi}^m, q)$ satisfies the following. 
\begin{equation}
\chi_{\quasi}(\scL_{\Phi}^m, q)=\chi_{\quasi}^{1/2}(\scL_{\Phi}^m, q)+
(-1)^\ell\chi_{\quasi}^{1/2}(\scL_{\Phi}^m, mh-q). 
\end{equation}
\end{proposition}
\begin{proof}

Write $R_{\Phi}^{1/2}(x)=\sum_{i=1}^{\lfloor h/2\rfloor}a_i'x^i$. 
From Theorem \ref{thm:yos-worp} and Proposition \ref{prop:dectrunc}, it 
follows that 
\begin{equation}
\chi_{\quasi}(\scL_{\Phi}^m, q)=
(R_{\Phi}^{1/2}(S^{m+1})L_\Phi)(q)+
(S^{h(m+1)}R_{\Phi}^{1/2}(S^{-m-1})L_\Phi)(q). 
\end{equation}
The first term is equal to $\chi_{\quasi}^{1/2}(\scL_{\Phi}^m, q)$. 
We shall prove that the second term is equal to 
$(-1)^\ell\chi_{\quasi}^{1/2}(\scL_{\Phi}^m, mh-q)$. Indeed, 
using Proposition \ref{prop:hshift}, we have 
\begin{equation}
\begin{split}
(S^{h(m+1)}R_{\Phi}^{1/2}(S^{-m-1})L_\Phi)(q)
&=
\sum_{i=1}^{\lfloor h/2\rfloor}a_i'S^{h(m+1)-(m+1)i}L_\Phi(q)\\
&=
\sum_{i=1}^{\lfloor h/2\rfloor}a_i'L_\Phi(q-h(m+1)+(m+1)i)\\
&=
(-1)^\ell \sum_{i=1}^{\lfloor h/2\rfloor}a_i'
L_\Phi(-q+hm-(m+1)i)\\
&=
(-1)^\ell \sum_{i=1}^{\lfloor h/2\rfloor}a_i'
S^{(m+1)i}L_\Phi(-q+hm)\\
&=
(-1)^\ell\chi_{\quasi}^{1/2}(\scL_{\Phi}^m, mh-q). 
\end{split}
\end{equation}
\end{proof}
The next corollary follows immediately. 

\begin{corollary}
$\chi_{\quasi}(\scL_{\Phi}^m, q)=(-1)^\ell
\chi_{\quasi}(\scL_{\Phi}^m, mh-q)$. 
\end{corollary}

Consider 
$F_\Phi(t):=
R_{\Phi}^{1/2}(S)t^\ell=\sum_{i=1}^{\lfloor h/2\rfloor}a_i'(t-i)^\ell$ 
as a polynomial in $t$. The distribution of the roots of $F_\Phi(t)=0$ 
will play a crucial role. 

\begin{proposition}
\label{prop:Remax}
Let $\Phi\in\{E_6, E_7, E_8, F_4, G_2\}$ be an exceptional root system. 
Suppose $F_\Phi(\alpha)=0, \alpha\in\C$. Then, $\Re\alpha<\frac{h}{2}$. 
\end{proposition}
\begin{proof}
This proposition can be verified computationally. 
We describe the method for the case $\Phi=E_6$. 
The other cases are similar. First, using Table \ref{fig:table}, 
Theorem \ref{thm:lp}, and 
$R_{A_6}(x)=x + 57 x^2 + 302 x^3 + 302 x^4 + 57 x^5 + x^6$, we have 
\[
\begin{split}
R_{E_6}(x)
=&
(1 + x)^3 (1 + x + x^2) 
(x + 57 x^2 + 302 x^3 + 302 x^4 + 57 x^5 + x^6)
\\
=&
x + 61 x^2 + 537 x^3 + 1916 x^4 + 3782 x^5 + 4686 x^6 + 3782 x^7 \\
&
+ 1916 x^8 + 537 x^9 + 61 x^{10} + x^{11}. 
\end{split}
\]
Therefore, we have 
\[
R_{E_6}^{1/2}(x)=x + 61 x^2 + 537 x^3 + 1916 x^4 + 3782 x^5 + 2343 x^6 
\]
and 
\[
R_{E_6}^{1/2}(S)t^6=(t - 1)^6 + 61 (t - 2)^6 + 537 (t - 3)^6 + 
1916 (t - 4)^6 +  3782 (t - 5)^6 + 2343 (t - 6)^6. 
\]
The roots of $R_{E_6}^{1/2}(S)t^6=0$ (approximation by Mathematica) are 
$
t=4.55334\pm 0.465487\sqrt{-1}, 4.78675\pm 1.55735\sqrt{-1}$, and 
$5.37033\pm 3.11072\sqrt{-1}$. All roots have real parts that are less than 
$\frac{h}{2}=6$. 

The maximum real parts of the roots are presented in Table \ref{fig:maxreal}. 
\begin{table}[htbp]
\centering
\begin{tabular}{c|c|c}
$\Phi$&$\max\{\Re\alpha\}$&$h/2$\\
\hline\hline
$E_6$&$5.3703$&$6$\\
$E_7$&$8.4367$&$9$\\
$E_8$&$14.6604$&$15$\\
$F_4$&$4.8967$&$6$\\
$G_2$&$2.166$&$3$
\end{tabular}
\caption{The maximal real parts of roots.}
\label{fig:maxreal}
\end{table}
\end{proof}

\section{A Toy Case}
\label{sec:toy}

Let $\Phi\in\{E_6, E_7, E_8, F_4, G_2\}$. Denote its exponents 
by $e_1, \dots, e_\ell$ and the Coxeter number by $h$. 
Let $g(t)\in\R[t]$ be a real polynomial of $\deg g=\ell$ that satisfies 
$g(t-h)=(-1)^\ell g(-t)$. (Note that such a polynomial exists, e. g., 
$g(t)=\prod_{i=1}^\ell(t+e_i)$.) 

\begin{theorem}
\label{thm:toyRH}
With the above notation, for sufficiently large $m\gg 0$, 
every root $\alpha\in\C$ of the equation $(R_\Phi(S^{m+1})g)(t)=0$ 
satisfies $\Re \alpha=\frac{mh}{2}$. 
\end{theorem}

\begin{proof}
Set $g_m(t):=(R_{\Phi}^{1/2}(S^{m+1})g)(t)$. Then, we have 
\[
(R_\Phi(S^{m+1})g)(t)=g_m(t)+(-1)^\ell g_m(mh-t). 
\]
(The proof is similar to that of Proposition \ref{prop:decquasi}.) 
By Proposition \ref{prop:normalim}, 
\[
\lim_{m\rightarrow\infty}
\frac{g_m(mt)}{m^\ell}=
R_{\Phi}^{1/2}(S)t^\ell=:F_\Phi(t). 
\]
From Proposition \ref{prop:Remax}, it follows that the real parts of 
the roots of the equation $F_\Phi(t)=0$ are less than $\frac{h}{2}$. 
By the continuity of the roots of polynomials, for sufficiently large 
$m\gg 0$, the real parts of the roots of $g_m(mt)=0$ are also less than 
$\frac{h}{2}$, which is equivalent to saying that the real parts of 
the roots of $g_m(t)=0$ are less than $\frac{mh}{2}$. Then 
Lemma \ref{lem:elem} completes the proof. 
\end{proof}

\section{Main Results}

\label{sec:main}

In this section, by generalizing the argument in \S \ref{sec:toy}, 
we will prove the main result. The main difficulty is related to 
the fact that $L_\Phi(q)$ is a quasi-polynomial. We will use 
the idea of averaging constituents to resolve this problem. 
We will work in the generalized setting 
described in \S \ref{subsec:settings} for the sake of notational simplicity. 

\subsection{Settings}
\label{subsec:settings}

Let $\widetilde{n}, h>0$ be positive integers and $L(q)$ be a 
quasi-polynomial with period $\widetilde{n}$. 

\begin{assumption}
\label{assump:1}
$L(q)$ has a constant leading term of degree $\ell>0$. In other words, 
$L(q)$ has an expression of the form 
\[
L(q)=c_0q^\ell+c_1(q)q^{\ell-1}+\cdots+c_\ell(q), 
\]
where $c_i:\Z\longrightarrow\Q$ ($i=1, \dots, \ell$) 
is a periodic function and $c_0\neq 0$. 
\end{assumption}

\begin{assumption}
\label{assump:2}
$L(q)$ satisfies the following. 
\begin{equation}
\label{eq:Ldual}
L(-q)=(-1)^\ell L(q-h). 
\end{equation}
\end{assumption}
Let $R(x)\in\Q[x]$ be a polynomial of $\deg R(x)=h-1$. 

\begin{assumption}
\label{assump:3}
$R(x)$ satisfies the following. 
\begin{equation}
\label{eq:Rdual}
x^hR(x^{-1})=R(x). 
\end{equation}
\end{assumption}
Write $R(x)=a_1x+a_2x^2+\cdots+a_{h-1}x^{h-1}$. 
Define the truncation $R'(x)$ of $R(x)$ by 
\[
R'(x)=
\left\{
\begin{array}{ll}
\sum\limits_{1\leq i<\frac{h}{2}}a_ix^i, & \mbox{ if $h$ is odd},\\
&\\
\sum\limits_{1\leq i<\frac{h}{2}}a_ix^i+\frac{a_{h/2}}{2}x^{h/2}, & \mbox{ if $h$ is even}. 
\end{array}
\right.
\]
We also write $R'(x)=\sum_{i=1}^{\lfloor h/2\rfloor}a_i'x^i$. It is 
easy to see that $R(x)$ satisfies 
\[
R(x)=R'(x)+x^hR'(x^{-1}). 
\]
Next, we make an assumption on the location of the roots of the 
polynomial 
\begin{equation}
R'(S)t^\ell=\sum_{i=1}^{\lfloor h/2\rfloor}a_i'(t-i)^\ell. 
\end{equation}

\begin{assumption}
\label{assump:4}
Every root $\alpha\in\C$ of $R'(S)t^\ell=0$ satisfies 
\begin{equation}
\Re\alpha<\frac{h}{2}. 
\end{equation}
\end{assumption}

The following is our main example. 
\begin{example}
\label{ex:mainex}
Let $\Phi\in\{E_6, E_7, E_8, F_4, G_2\}$. Then, the Ehrhart quasi-polynomial 
$L_\Phi(q)$, period $\widetilde{n}$, Coxeter number $h$, and 
Eulerian polynomial $R_\Phi(x)$ satisfy 
Assumptions 
\ref{assump:1}, 
\ref{assump:2}, 
\ref{assump:3}, and 
\ref{assump:4}. 
\end{example}

\begin{remark}
Although Example \ref{ex:mainex} is the main example, 
we can construct many other examples that satisfy the above assumptions. 
For instance, for $\widetilde{n}=1$, some arbitrary $h\geq 3$ and $\ell>0$, 
$L(q)=(q+\frac{h}{2})^\ell$ and $R(x)=x+x^{h-1}$ satisfy 
the above assumptions. 
\end{remark}

\subsection{Asymptotic behavior of roots}

\label{subsec:asympt}

For $m>0$, define the quasi-polynomials $L^{(m)}(q)$ and $L'^{(m)}(q)$ 
of period 
$\widetilde{n}$ by 
\[
L^{(m)}(q)=(R(S^{m+1})L)(q) 
\]
and 
\[
L'^{(m)}(q)=(R'(S^{m+1})L)(q). 
\]
By Assumption \ref{assump:2} (and Assumption \ref{assump:3}), 
\begin{equation}
\label{eq:decLprim}
L^{(m)}(q)=L'^{(m)}(q)+(-1)^\ell L'^{(m)}(mh-q). 
\end{equation}

Let us denote the constituents by $L_d^{(m)}(t)\in\Q[t]$ 
for each residue class $d\mod\widetilde{n}$ 
(or $0\leq d<\widetilde{n}$), namely, 
\[
L^{(m)}(q)=
L_d^{(m)}(q), \mbox{ when }q\equiv d\mod\widetilde{n}. 
\]
Note that the constituents have the following expression: 
\[
L_d^{(m)}(t)=
\sum_{i=1}^{h-1}
a_i\sum_{j=0}^\ell c_j(d-(m+1)i)\cdot (t-(m+1)i)^{\ell-j}. 
\]
The relation (\ref{eq:decLprim}) can also be written as 
\begin{equation}
\label{eq:decconst}
L_d^{(m)}(t)=L_d'^{(m)}(t)+(-1)^\ell\cdot L_{mh-d}'^{(m)}(mh-t) 
\end{equation}
at the level of the constituents. 

\begin{proposition}
\label{prop:Lprim}
\[
\lim_{m\rightarrow\infty}\frac{L_d'^{(m)}(mt)}{m^\ell}=c_0\cdot R'(S)t^\ell. 
\]
In particular, the limit does not depend on the residue $d$. 
\end{proposition}

\begin{proof}
Recall $L_d'^{(m)}(t)=
\sum_{i=1}^{\lfloor h/2\rfloor}a_i'
\sum_{j=0}^\ell c_j(d-(m+1)i)(t-(m+1)i)^{\ell-j}$. Hence, 
\[
\frac{L_d'^{(m)}(mt)}{m^\ell}=
\sum_{i=1}^{\lfloor h/2\rfloor}a_i'
\sum_{j=0}^\ell \frac{c_j(d-(m+1)i)}{m^j}
\left(t-\frac{m+1}{m}i\right)^{\ell-j}. 
\]
When $j>0$, since $c_j$ is a periodic function, we have 
$\lim_{m\rightarrow\infty}\frac{c_j(d-(m+1)i)}{m^j}=0$. 
By the assumption, $c_0(d-(m+1)i)= c_0$ is a nonzero constant. 
We have 
\[
\begin{split}
\lim_{m\rightarrow\infty}\frac{L_d'^{(m)}(mt)}{m^\ell}
&=
c_0\cdot\sum_{i=1}^{\lfloor h/2\rfloor}a'_i(t-i)^\ell\\
&=
c_0\cdot R'(S)t^\ell. 
\end{split}
\]
\end{proof}

\begin{definition}
\label{def:infsup}
\[
\begin{split}
\overline{r}_d^{(m)}&:=
\max\{\Re\alpha\mid \alpha\in\C, L_d^{(m)}(\alpha)=0\}, \\
\underline{r}_d^{(m)}&:=
\min\{\Re\alpha\mid \alpha\in\C, L_d^{(m)}(\alpha)=0\}. 
\end{split}
\]
\end{definition}
If $L(q)$ and $R(x)$ satisfy Assumptions 
\ref{assump:1}--\ref{assump:4}, then 
$\overline{r}_d^{(m)}$ and 
$\underline{r}_d^{(m)}$ approach $\frac{mh}{2}$ as $m\rightarrow\infty$. 
More precisely, we have the following. 

\begin{theorem}
\label{thm:asymptr}
For any $0\leq d<\widetilde{n}$, 
\begin{equation}
\label{eq:limithhalf}
\lim_{m\rightarrow\infty}
\frac{\overline{r}_d^{(m)}}{m}=
\lim_{m\rightarrow\infty}
\frac{\underline{r}_d^{(m)}}{m}=\frac{h}{2}. 
\end{equation}
\end{theorem}

\begin{proof}
Let $F(t):=c_0\cdot R'(S)t^\ell$. Then, by (\ref{eq:decconst}) and 
Proposition \ref{prop:Lprim}, we have 
\begin{equation}
\lim_{m\rightarrow\infty}
\frac{L_d^{(m)}(mt)}{m^\ell}=F(t)+(-1)^\ell\cdot F(h-t). 
\end{equation}
Choose a root $\alpha_m\in\C$ of $L_d^{(m)}(t)=0$. 
Then, obviously, $\frac{\alpha_m}{m}$ satisfies the equation 
$L_d^{(m)}(mt)=0$. Hence, $\frac{\alpha_m}{m}$ approaches 
the (set of) roots of $F(t)+(-1)^\ell F(h-t)=0$. 
Equation (\ref{eq:limithhalf}) follows from 
Assumption \ref{assump:4} and Lemma \ref{lem:elem}. 
\end{proof}

\begin{corollary}
\label{cor:limRH}
Let $\Phi\in\{E_6, E_7, E_8, F_4\}$. Fix $0\leq d<\widetilde{n}$. 
Let $\alpha_m\in\C$ be a root of the constituent 
$\chi_{\quasi, d}(\scL_\Phi^m, t)$. Then, 
\[
\lim_{m\rightarrow\infty}
\frac{\Re\alpha_m}{m}=\frac{h}{2}. 
\]
\end{corollary}

\subsection{Exact arrangement of roots}

\label{subsec:exact}

In the previous subsection, we proved that the real part of a root 
$\alpha$ of a constituent $L_d^{(m)}(t)=0$ is asymptotically close to 
$\frac{mh}{2}$. 
Here, we prove a stronger result for some special constituents. 

\begin{definition}
\label{def:adm}
(Using the notation in \S \ref{subsec:settings}), 
the residue $d\mod \widetilde{n}$ ($0\leq d<\widetilde{n}$) is called 
an \emph{admissible residue} if the constituents of $L^{(m)}(q)$ satisfy 
\begin{equation}
\label{eq:adm}
L_d^{(m)}(t)=
L_{d+kh}^{(m)}(t)=
L_{-d+kh}^{(m)}(t) 
\end{equation}
for all $k\in\Z$ and $m>0$. 
\end{definition}

\begin{example}
Consider the situation in Example \ref{ex:mainex}. Then, 
$L^{(m)}(q)=(R(S^{m+1})L)(q)$ is equal to $\chi_{\quasi}(\scL_\Phi^m, q)$, 
which has the $GCD$-property. Hence, $d$ is an admissible residue 
if and only if 
\[
\gcd(d, \widetilde{n})=
\gcd(d+kh, \widetilde{n}) 
\]
for all $k\in\Z$. As $\rad(\widetilde{n})|h$ (Theorem \ref{thm:suter} 
(v)), $d=1$ is an admissible residue. Other admissible residues (divisors 
of $\widetilde{n}$) for $\Phi\in\{E_6, E_7, E_8, F_4, G_2\}$ 
are listed in Table \ref{tab:ar}. 
\begin{table}[htbp]
\centering
\begin{tabular}{c|c|c|c|c|c}
$\Phi$&$\widetilde{n}$&$\rad(\widetilde{n})$&$h$&\mbox{admissible divisor 
of $\widetilde{n}$}&$m_0$\\
\hline\hline
$E_6$&$6$&$6$&$12$&$1,2,3,6$&$1$\\
$E_7$&$12$&$6$&$18$&$1,3$&$2$\\
$E_8$&$60$&$30$&$30$&$1,3,5,15$&$2$\\
$F_4$&$12$&$6$&$12$&$1,2,3,4,6,12$&$1$\\
$G_2$&$6$&$6$&$6$&$1,2,3,6$&$1$
\end{tabular}
\caption{Admissible divisors.}
\label{tab:ar}
\end{table}
\end{example}

The following is the main result. 

\begin{theorem}
\label{thm:maingeneral}
(Using the same notation as in \S \ref{subsec:settings}.) 
Suppose $d$ ($0\leq d<\widetilde{n}$) is an admissible residue. 
Let $\alpha_m$ be a root of $L_d^{(m)}(t)=0$. Then, for sufficiently 
large $m\gg 0$, $\Re\alpha_m=\frac{mh}{2}$ holds. 
\end{theorem}

\begin{proof}
Let $m_0:=\frac{\widetilde{n}}{\gcd(h, \widetilde{n})}$. Then, 
$L_d^{(m)}(t)=L_{d+kh}^{(m)}(t)=L_{-d+kh}^{(m)}(t)$ for 
$k=0, 1, \dots, m_0-1$. Hence, 
using (\ref{eq:decconst}), 
\begin{equation}
\label{eq:2m01}
\begin{split}
L_d^{(m)}(t)
&=
\frac{1}{m_0}\cdot\sum_{k=0}^{m_0-1}L_{d+kh}^{(m)}(t)\\
&=
\frac{1}{m_0}\cdot\sum_{k=0}^{m_0-1}L_{d+kh}'^{(m)}(t)
+\frac{(-1)^\ell}{m_0}\cdot\sum_{k=0}^{m_0-1}L_{-d+(m-k)h}'^{(m)}(mh-t)\\
&=
\frac{1}{m_0}\cdot\sum_{k=0}^{m_0-1}L_{d+kh}'^{(m)}(t)
+\frac{(-1)^\ell}{m_0}\cdot\sum_{k=0}^{m_0-1}L_{-d+kh}'^{(m)}(mh-t). 
\end{split}
\end{equation}
As $L_d^{(m)}(t)=L_{-d}^{(m)}(t)$, we also have 
\begin{equation}
\label{eq:2m02}
L_d^{(m)}(t)
=
\frac{1}{m_0}\cdot\sum_{k=0}^{m_0-1}L_{-d+kh}'^{(m)}(t)
+\frac{(-1)^\ell}{m_0}\cdot\sum_{k=0}^{m_0-1}L_{d+kh}'^{(m)}(mh-t). 
\end{equation}
Define the polynomial $F_d^{(m)}(t)$ by 
\[
F_d^{(m)}(t):=
\frac{1}{2m_0}\cdot
\sum_{k=0}^{m_0-1}
\left\{
L_{d+kh}'^{(m)}(t)+L_{-d+kh}'^{(m)}(t)
\right\}. 
\]
Then, combining (\ref{eq:2m01}) and (\ref{eq:2m02}), 
$L_d^{(m)}(t)$ can be expressed as 
\[
L_d^{(m)}(t)=F_d^{(m)}(t)+(-1)^\ell\cdot F_d^{(m)}(mh-t). 
\]
From Proposition \ref{prop:Lprim}, it follows that 
\begin{equation}
\lim_{m\rightarrow\infty}
\frac{F_d^{(m)}(mt)}{m^\ell}=c_0\cdot R'(S)t^\ell. 
\end{equation}
Hence, for sufficiently large $m\gg 0$, every root $\alpha\in\C$ of 
$F_d^{(m)}(t)=0$ satisfies 
\[
\Re\alpha<\frac{mh}{2}. 
\]
Applying Lemma \ref{lem:elem}, every root of 
$L_d^{(m)}(t)=0$ has the real part $\frac{mh}{2}$. 
\end{proof}

\begin{corollary}
\label{cor:main}
Let $\Phi\in\{E_6, E_7, E_8, F_4, G_2\}$. For sufficiently large 
$m\gg 0$, every root $\alpha\in\C$ of the characteristic polynomial 
$\chi(\scL_\Phi^m, t)$ of the Linial arrangement $\scL_\Phi^m$ 
satisfies $\Re\alpha=\frac{mh}{2}$. 
\end{corollary}

\begin{proof}
Recall that the characteristic polynomial is the constituent 
of the characteristic quasi-polynomial corresponding to $d=1$. 
Since $d=1$ is an admissible residue, 
we can apply Theorem \ref{thm:maingeneral}. 
\end{proof}


\medskip

\noindent
{\bf Acknowledgements.} 
The author was partially supported by 
JSPS KAKENHI Grant Number 
25400060, 
15KK0144, 
and 
16K13741. 

\end{document}